\documentclass[11pt,a4paper]{article}

\usepackage{amsmath} 
\usepackage{amssymb}
\usepackage{dsfont}
\usepackage{graphicx}
\usepackage{here}
\newtheorem{theorem}{Theorem}[section]

\newtheorem{definition}[theorem]{Definition}
\newtheorem{remark}{Remark}

\newcommand{\D}{{\, \rm d}}

\newcommand{\R}{\mathbb{R}}

\newcommand{\bu}{\mathbf{u}}
\newcommand{\bx}{\mathbf{x}}
\newcommand{\be}{\mathbf{e}}
\newcommand{\mH}{\mathcal{H}}

\usepackage{xspace,xcolor}
\usepackage[linecolor=blue!40!white, backgroundcolor=blue!05!white,bordercolor=blue!40!white,textsize=tiny]{todonotes}
\let\origintodo\todo  \newcommand{\xtodo}[2][]{\origintodo[#1]{#2}\xspace}  \let\todo\xtodo

\title{Numerical simulation on a cell polarisation model: the polar case}

\author{Vincent Calvez\thanks{Unit\'e de Math\'ematiques Pures et Appliqu\'ees, CNRS UMR 5669 \& \'equipe-projet INRIA NUMED, \'Ecole Normale Sup\'erieure de Lyon, 46 all\'ee d'Italie, F-69364 Lyon, France.  ({\tt vincent.calvez@umpa.ens-lyon.fr})} \and Nicolas Meunier\thanks{MAP5, CNRS UMR 8145, Universit\'{e} Paris Descartes, 45 rue des Saints  P\`{e}res
75006 Paris,
France. ({\tt nicolas.meunier@parisdescartes.fr})} \and Nicolas Muller \thanks{MAP5, CNRS UMR 8145, Universit\'{e} Paris Descartes, 45 rue des Saints  P\`{e}res
75006 Paris,
France. {\tt nicolas.muller@parisdescartes.fr})}  \and Raphael Voituriez.
        \thanks{Laboratoire de la mati\`ere condens\'ee, CNRS UMR 7600, Universit\'e Pierre et Marie Curie, 4 Place Jussieu, 75255 Paris Cedex 05 France ({\tt voiturie@lptmc.jussieu.fr})}}

\begin{document}

\maketitle

\begin{abstract}
When it is polarised, a cell develops an asymmetric distribution of specific molecular markers, cytoskeleton and cell membrane shape. Polarisation can occur spontaneously or be triggered by external signals, like gradients of signalling molecules... In this work, we use the models of cell polarisation introduced in \cite{Firstpaper} and \cite{Siam_CHMV} and we set a numerical analysis for these models. They are based on nonlinear convection-diffusion equations and the nonlinearity in the transport term expresses the positive loop between the level of protein concentration localised in a small area of the cell membrane and the number of new proteins that will be convected to the same area. We perform numerical simulations and we illustrate that these models are rich enough to describe the apparition of a polarisome.   
\end{abstract}

\section{Introduction}
Cell polarisation is a symmetry-breaking event that occurs in cell division, mating or morphogenesis. Molecular markers play a central role in establishing this phenomenon. Indeed, there are two different behaviours: a non-polarised cell has its markers radially homogeneously distributed while markers are located in a small area of the cell membrane for a polarised cell. Yeast cells are dynamically polarised in response to the extracellular gradients of signals (chemokynes). However, it has been observed in \cite{Altschuler2003} that polarisation can occur spontaneously without any external asymmetric stimulus. 

During the past decade, many models describing cell polarisation have been developed. The majority of these models are based on reaction-diffusion systems where polarisation appears as a type of Turing instability \cite{Turing1}, \cite{Turing2}, \cite{Turing3}, or due to stochastic fluctuations \cite{Altschuler}, other models include cytoskeleton proteins as a regulatory factor \cite{Marco2007}, \cite{Altschuler2003}. Many biological studies have shown that the cytoskeleton plays an important role in polarisation. It has been suggested that there is a positive feedback on molecular markers density. Indeed, disruption of transport along the cytoskeleton greatly reduces the stability of polar cap \cite{Altschuler2003}. The cell cytoskeleton is a network of long semi-flexible filaments made up of protein subunits \cite{cytoskeleton}. These filaments (mainly actin or microtubules) act as roads along which motor proteins are able to perform a biased ballistic motion and carry various molecules. Molecular markers play a key role in the formation of these filaments.

Following \cite{Firstpaper}, \cite{CRAS} and \cite{Siam_CHMV}, in this work we study models that describe the dynamics of cell polarisation. In these models, molecular markers, such as proteins, diffuse in the cytoplasm and are actively transported along the cytoskeleton. The resulting motion is a biased diffusion regulated by the markers themselves. Using numerical simulations and mathematical heuristics, we observe that the coupling on the velocity field achieves an inhomogeneous distribution of molecular markers without any external asymmetric field. Such an inhomogeneous distribution is only due to interaction between molecular markers.

Throughout this paper, the density of molecular markers (resp. advection field) is denoted by $\rho(t,\bx)$ (resp. $\bu(t,\bx)$). The advection is obtained through a coupling with the membrane concentration of markers. The cell is figured by the domain $\Omega \subset \mathbb{R}^n$ with $n=1,2$ and a part of the boundary of the domain will be the active membrane denoted by $\Gamma$.  The time evolution of the molecular markers satisfies the following advection-diffusion equation, see \cite{Firstpaper} and \cite{Siam_CHMV}:
\begin{equation}\label{main}
\left\{
\begin{aligned}
\partial_t \rho(t,\bx) & = D\, \Delta \rho(t,\bx) - \chi\, \nabla . \left( \rho(t,\bx) \, \bu (t,\bx) \right), \quad t > 0, \quad \bx \in \Omega, \\
\rho(0,\bx) & = \rho_0 (\bx).
\end{aligned}
\right.
\end{equation}
There is no creation nor degradation of molecular markers in the cell, so the quantity of molecular markers remains constant in time:
\begin{equation}\label{massconservation}
M = \int_{x\in\Omega} \rho_0(\bx) d\bx = \int_{x\in\Omega} \rho(t,\bx) d\bx.
\end{equation}
This condition is ensured by a zero flux boundary condition on the boundary.
A first simplified step is to assume that the cell is essentially bidimensional and to neglect curvature effects. The membrane
boundary is then a 1D line along the $y$-axis and the cytoplasm is parametrized by $\bx=(x,y) \in \mathbb{R}_+ \times \mathbb{R}$.

The plan of this work is the following. First, we present the models, that are based on different expressions for $\bu$, and we recall the main mathematical results of the simplified model in 1D for $\Omega=(0,\infty)$ and $\Gamma = \{ x = 0 \}$, see \cite{CRAS}, \cite{Siam_CHMV} for more details. Then we study a more realistic model, that includes dynamical exchange of markers on the boundary. This model was introduced in \cite{Firstpaper} and studied in \cite{Siam_CHMV}. We provide a methodology for parameter estimation and qualitative description of cell polarisation  by using mathematical heuristics. Then, we perform a numerical analysis of this  model. We introduce the numerical part by the one dimensional case. Finally, we give tools to study the numerical implementation of the model on an annulus domain.

\section{Presentation of the models and mathematical results}

\subsection{One dimensional case}
In this section, we study the one dimensional case on the half line for $\Omega=(0,\infty)$. The membrane is then the point $\Gamma = \{ x = 0 \}$. For the first model, the advection field towards the membrane is equal to the density of molecular markers on the boundary $\rho(t,0)$. Then we improve this model by considering that only the trapped molecular markers on the membrane contribute to the advection field.
\subsubsection{Simplified model set on the half line}
In \cite{Siam_CHMV} a first mathematical studies  has been done on this model. We define an advection field $\bu(t,x)$ for \eqref{main} $$\bu(t,x)=-\rho(t,0),$$
in such a case \eqref{main} reads as (with $D=1$ and $\chi=1$):
\begin{equation}\label{simplified}
\partial_t \rho(t,x) = \partial_{xx} \rho(t,x) + \rho(t,0) \, \partial_x \rho(t,x), \quad t > 0, \quad x>0,
\end{equation}
with the following zero flux condition on the boundary $\Gamma=\{x=0\}$, that ensures the mass conversation \eqref{massconservation},
\begin{equation}
\partial_x \rho(t,0) + \rho(t,0)^2=0.
\end{equation}
Interestingly enough in \cite{Siam_CHMV}, it has been proved that solutions of \eqref{simplified} blow-up in finite time if their masses are above a certain critical mass, $ M> 1 $, and exist globally in time if $ M \leq 1 $. Let us first recall the definition of weak solutions of \eqref{simplified}.
\begin{definition}\label{def:weak}
We say that $\rho(t,x)$ is a weak solution of \eqref{simplified}  on $(0,T)$ if it satisfies:
\begin{equation*}
\rho \in L^\infty(0,T;L^1_+(\R_+))\, , \quad \partial_x \rho \in L^1((0,T)\times \R_+)  \, , \label{eq:flux L1}
\end{equation*}
and $\rho(t,x)$ is a solution of \eqref{simplified} in the sense of distributions in $\mathcal D'(\R_+)$.
\end{definition}

Let us now recall the main results for weak solutions of \eqref{simplified}.
\begin{theorem}[Global existence: $M\leq1$] \label{th:1D} 
Assume that the initial data $\rho_0$ satisfies both $\rho_0 \in L^1(( 1 + x)\D x)$ and $\int_{x>0} \rho_0(x) (\log \rho_0(x))_+ \D x< + \infty$. Assume in addition that $M\leq 1$, then there exists a global weak solution of equation \eqref{simplified}.
\end{theorem}
\begin{theorem}[Blow-up: $M>1$] \label{th:1D BU} Assume $M>1$. Any weak solution of equation \eqref{simplified} with non-increasing initial data $\rho_0$ blows-up in finite time.
\end{theorem}
\begin{remark}
It would tempted to interpret blow-up of solutions of the one dimensional model as cell polarisation but concentration of markers on the boundary doesn't automatically mean polarisation. Indeed, consider a radially symmetric 2D cell case then equation \eqref{main} reduces to the one dimensional one. Above a threshold on the total mass, the convection wins and markers concentrate on the boundary. In some situations, these markers may be homogeneously distributed on the boundary and in such a case there is no symmetry breaking.
\end{remark}

\subsubsection{The model with dynamical exchange of markers at the boundary} 
Such a direct activation of transport on the boundary seems to be unrealistic. Indeed possible occurrence of blow-up in finite time suggests this claim. We improve the previous model by distinguishing between cytoplasmic content $\rho(t,x)$ and the concentration of trapped molecules on the boundary, that will be denoted by $\mu(t)$. The dynamical exchange of markers at the boundary is done with an attachment rate $k_{on}$ and a detachment rate $k_{off}$, hence the time evolution of $\mu(t)$ is
\begin{equation}\label{mu1D}
\dfrac{d}{dt} \mu (t) = k_{on} \, \rho(t,0) - k_{off} \, \mu(t).
\end{equation}
The advection field $\bu(t,x)$ in \eqref{main} is now defined by $$\bu(t,x)=-\mu(t),$$
hence \eqref{main}  (with $D=1$ and $\chi=1$) reads as:
\begin{equation}
\partial_t \rho(t,x) = \partial_{xx} \rho(t,x) + \mu(t) \, \partial_x \rho(t,x), \quad t > 0, \quad x>0,
\end{equation}
with a modified boundary condition
\begin{equation}
\partial_x \rho(t,0) + \rho(t,0)\, \mu(t) = \dfrac{d}{dt} \mu (t).
\end{equation}
This ensures the following mass conservation shared among $\rho(t,x)$ and $\mu(t)$:
\begin{equation}
M = \int_{\mathbb{R}_+} \rho_0(x) dx + \mu_0 = \int_{\mathbb{R}_+} \rho(t,x) dx + \mu(t).
\end{equation}
With equation \eqref{mu1D}, the self-activation of transport by $\rho(t,0)$ is then delayed in time. Since the transport speed is bounded $\mu(t) \leq M$, the solution of the model with dynamical exchange on the boundary exists globally in time. More precisely it is possible, see \cite{Siam_CHMV}, to prove that it converges towards a non trivial stationary state. 
\begin{theorem}[Global existence: dynamical exchange case] \label{th:1Dmu} 
Assume that the initial data $\rho_0$ satisfies both $\rho_0 \in L^1(( 1 + x)\D x)$ and $\int_{x>0} \rho_0(x) (\log \rho_0(x))_+ \D x< + \infty$. Assume the mass is super-critical $M > 1$.
The partial mass $m(t)=\int_{x>0} \rho(t,x) \D x$ converges to $1$ and the density $\rho(t,x)$ strongly converges in $L^1$ towards the exponential profile $(M -1)\exp(-(M-1)x)$.
\end{theorem}

\subsection{Two dimensional case} 
In this section, we study the two dimensional case on $\Omega  \subset \mathbb{R}^2$. As in the one dimensional case, the advection field towards the membrane depends on the density of molecular markers on the boundary, two different situations (actin and microtubules) will be described,  then we improve this model by considering an exchange of markers at the membrane and only the trapped molecular markers contribute to the advection field.
\subsubsection{Simplified model set on the half plan}
We study the model on the half plan for $\bx = (x,y) \in \Omega=\mathbb{R}_+ \times \mathbb{R}$. The membrane is then the line $\Gamma = \{ x = 0 \} \times \mathbb{R}$. We have the following boundary condition for $\rho(t,\bx)$ at point $\bx \in \Gamma$
\begin{equation}
(D\, \nabla \rho(t, \bx) - \chi \, \rho(t,\bx)\, \bu(t,\bx) ). \vec n_\bx  = 0,
\end{equation}
where $\vec n_\bx$ is the outward normal to $\Gamma$.
This ensures the following mass conservation:
\begin{equation}
M = \int_{\Omega} \rho_0(\bx) d\bx = \int_{\Omega} \rho(t,\bx) d\bx.
\end{equation}
First we consider the transversal case, the field $\bu$ is normal to the boundary
\begin{equation}\label{MT}
\bu(t,\bx) = - S(y) \rho(t,0,y) \vec \be_x.
\end{equation}
The microtubules of the cytoskeleton are normally oriented to the cell membrane and their growth depends on the density of molecular markers on the boundary. In the potential case, we consider the following advection field deriving from a harmonic potential modelling the transport by actin filaments: 
\begin{equation}\label{c}
\bu(t,\bx) = \nabla c(t,\bx), \mbox{ where } \begin{cases}- \Delta c(t,\bx) = 0, & \mbox{ if } \bx \in \Omega, \\
\nabla c(t,\bx) . \vec n_\bx = S(\bx) \rho(t,\bx), & \mbox{ if } \bx \in \Gamma.
\end{cases}
\end{equation}
This advection field orientation is due to the actin networks.
\begin{center}
\includegraphics[scale=0.7]{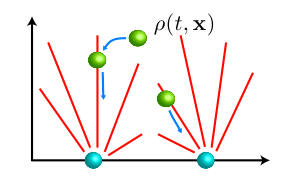}
\end{center}
Actin filaments are attached on the membrane and randomly distributed, there orientations are mixed up. We also add the external pheromone concentration at $\bx \in \Gamma$ which acts by the mating-pheromone MAPK cascade on the actin transport.
We have global existence and blow-up theorems for the simplified model with $\bx=(x,y) \in \Omega = \R_+ \times \R$. For clarity, we recall this result, see \cite{Siam_CHMV} for more details. 
\begin{theorem}[Global existence in dimension 2]\label{thdim2}
Assume that the advection field satisfies the two following conditions: $\nabla\cdot \bu \geq 0$ and $\bu(t,0,y)\cdot \vec \be_x = - \rho(t,0,y)$. Assume that the initial data $\rho_0$ satisfies both $\rho_0 \in L^1(( 1 + |\bx|^2)\D \bx)$ and $\|\rho_0\|_{L^2}$ is smaller than some constant $c$. Then there exists a global weak solution to equations \eqref{main}-\eqref{massconservation}. 
\end{theorem}
\begin{theorem}[Blow-up in dimension 2]\label{thdim2bu} Assume that $\rho(t, \bx)$ is a strong solution to \eqref{main} which verifies: 
\begin{itemize}
\item
$\partial_x \rho(t, \bx) \leq 0$ for all $\bx \in \Omega$  and $t > 0$ when the advective field is given by \eqref{MT}.
\item
$\partial_x \rho(t, \bx) \leq 0$ for all $\bx \in \Omega$  and $t > 0$ the matrix $A(t, \bx) = \bx \otimes \partial_x \partial_y \log \rho(t, \bx)$ satisfies $A^T +A \geq 0$ (in the matrix sense) when the advective field is given by \eqref{c}.
\end{itemize}
Assume in addition that the second momentum is initially small enough: there exists a constant $C$ such that $ \int_{\bx \in \Omega} |\bx| ^2 \rho_0(\bx) d\bx < C M^3$. Then the maximal time of existence of the solution is finite.
\end{theorem}

\subsubsection{The model with dynamical exchange of markers at the boundary} 
Let $\Omega \subset \mathbb{R}^2$ be the cytoplasm domain, as in the one dimensional case \eqref{mu1D} we consider dynamical exchange of markers at the boundary, so for $\bx \in \Gamma$ we have the evolution in time of $\mu(t,\bx)$
\begin{equation}
\partial_t \mu (t,\bx) = k_{on} \, \rho(t,\bx) - k_{off} \, \mu(t,\bx),
\end{equation}
with a modified boundary condition for $\rho(t,\bx)$ at point $\bx \in \Gamma$
\begin{equation}
(D\, \nabla \rho(t, \bx) - \chi \, \rho(t,\bx)\, \bu(t,\bx) ). \vec n_\bx  = - \partial_t \mu (t,\bx),
\end{equation}
where $\vec n_\bx$ is the outward normal to $\Gamma$.
This ensures the following mass conservation sharing by $\rho(t,\bx)$ and $\mu(t,\bx)$:
\begin{equation}
M = \int_{\Omega} \rho_0(\bx) d\bx + \int_{\Gamma} \mu_0(\bx) d\bx = \int_{\Omega} \rho(t,\bx) d\bx + \int_{\Gamma} \mu(t,\bx) d\bx.
\end{equation}
As before, we consider the following advection field (cytoskeleton): 
\begin{equation}
\bu(t,\bx) = \nabla c(t,\bx), \mbox{ where } \begin{cases}- \Delta c(t,\bx) = 0, & \mbox{ if } \bx \in \Omega, \\
\nabla c(t,\bx) . \vec n_\bx = S(\bx) \mu(t,\bx), & \mbox{ if } \bx \in \Gamma.
\end{cases}
\end{equation}
We can also consider the transversal case. For example for $\Omega = \mathbb{R}_+ \times \mathbb{R}$, we take similarly as before
\begin{equation}
\bu(t,\bx) = - S(y) \mu(t,y) \vec{\be_x}.
\end{equation}
For the model with exchange on the boundary,  blow-up or global existence have not been proved yet. In this work, we make a first step in this direction by using a mathematical heuristic and numerical simulations. There is also one open question: does advection field (transversal or potential) create a break of symmetry ?

\section{Heuristic}\label{Heuristic}

The mathematical analysis performed in \cite{Siam_CHMV} has demonstrated that a class of models exhibit pattern formation (either blow-up or convergence towards a non homogeneous steady state) under some conditions. However the main question still remains unanswered: do these models describe cell polarisation or not? Thus in order to provide a first answer to this question, in the next section, we will perform numerical simulations. Our aim will be to see if, under some conditions, the model leads to a concentration of markers, not only on the boundary, but on a small region of the boundary.

Let us now describe the mathematical heuristic which is done on the half plan $\Omega = \mathbb{R}_+ \times \mathbb{R}$. We give formal arguments to motivate the differences arising in the dynamics leading by the two possible drifts $\bu_T$ (transversal case) and $\bu_P$ (potential case). Let $\bu_P=\nabla c$ be the solution of \eqref{c} we have (see \cite{Evans})
$$c(x,y) = - \frac{1}{\pi} \int_{y'\in \mathbb{R}}\log (\sqrt{(y-y')^2+x^2}) S(y') \rho(t,0,y') dy'.$$
We notice that the two possible drifts $\bu_T$ and $\bu_P$  share common features: they are both divergence free and their normal components at the boundary coincide.
$$\bu_T \cdot \vec{\be_x} =  \bu_P \cdot \vec{\be_x} = - S(y) \rho(t,0,y) \, .$$
On the other hand, a key difference holds when looking at the tangential component at the boundary: 
\begin{align*}  &\bu_T \cdot \vec{\be_y} = 0 \, , \\ 
&\bu_P \cdot \vec{\be_y} = - \pi \mH (S(\cdot) \rho(t,0,\cdot))\, ,
\end{align*}
where $\mH$ denotes the one-dimensional Hilbert transform:
\begin{equation*}
\mH (f) (y) = \dfrac1\pi{\rm p.v.}\int_{\R}\dfrac1{y-x} f(x)\, dx\, .
\end{equation*}
We expect the solution to concentrate on the boundary in the super-critical case for both choices of $\bu$, numerical simulations suggest it (see section \ref{Graphics}). Postulating the ansatz $\rho(t,x,y) = \nu(t,y)\delta(x = 0)$, we can formally write the dynamics of $\nu(t,y)$ for the two cases. Integrating the main equation \eqref{main} with respect to $x$ with zero flux condition on $\Gamma = \{x=0\}$, we obtain:
\begin{equation*}
\partial_t \int_{\mathbb{R}_+} \rho(t,x,y) dx  = D \partial_{yy} \left(\int_{\mathbb{R}_+} \rho(t,x,y) dx\right) - \chi \partial_y \left(\int_{\mathbb{R}_+} \rho(t,x,y) (\bu(t,x,y) \cdot \vec \be_y) dx \right).
\end{equation*}
\begin{itemize}
\item $\bu_T \cdot \vec{\be_y} = 0 $: diffusion equation for $\nu$:
\begin{equation*}
\partial_t \nu(t,y)  = D \, \partial_{yy} \nu(t,y)\, .
\end{equation*} 
\item $\bu_P \cdot \vec{\be_y} = - \pi\mH (S(\cdot) \rho(t,0,\cdot)) $ and assuming $S$ constant on $\mathbb{R}$, it reads as:
\begin{equation*}
\partial_t \nu(t,y)  = D \, \partial_{yy} \nu(t,y) + \chi   S  \, \partial_y \left(\nu(t,y) \mH ( \nu) (y) \right)\, .
\end{equation*}
\end{itemize}
Transversal dynamics are very different: in transversal case boundary diffusion dominates while in potential case
Hilbert transform has a critical singularity to offset the diffusion on this equation. This latter equation exhibits blow-up if $\int_\R \nu(t,y)\, dy=M$ is above the critical mass $\frac{2 \pi D}{ S \chi }$, see \cite{KSHilbert} done in a peculiar variant of Keller-Segel equation for more informations.  
\begin{remark}
This is a first step to observe a critical mass phenomenon and this may lead to blow-up if the mass is large enough. In this way, we also define an order of magnitude for some parameters. 
It is to be noticed that this latter criterion is valid for an infinite domain, namely $y\in \R$. In the case of a cell, the domain is finite and the existence of such a dichotomy has not been proved yet. In order to obtain more information on the critical value distinguishing the polarised case and the stable case we will perform numerical simulations.
\end{remark}

This heuristic shows us that potential case is able to break symmetry more readily than the transversal case. This difference between the two models is also discussed in the following section.

\section{Numerical analysis}

We introduce this numerical section by the discretization of the convection-diffusion model set on a 1D periodic domain. This first step allows us introducing the discretization of this model on a 2D domain. We model the cell as an annulus, molecular markers cannot pass inside the nucleus. In this section, for simplicity we fix all the parameters values to 1 except $M$.

\subsection{First step with the one dimensional case} Let $u(t,\theta)$ be a periodic function on the periodic domain $\Omega = \mathbb{R}/2\pi\mathbb{Z}$. We consider the following advection-diffusion equation 
\begin{equation}\label{1Dpolar}
\partial_t \rho(t,\theta) =  \partial_{\theta}\left (\partial_{\theta}\rho(t,\theta) - \rho(t,\theta) u(t,\theta)\right), \quad t > 0, \quad \theta\in \mathbb{R}/2\pi\mathbb{Z}.
\end{equation}
Let $t^n = n\, \Delta t$ be the time discretization and $\{\theta_k=k \, \Delta \theta, k \in \{1,...,N_\theta\}\}$ be the space discretization of the periodic interval $\mathbb{R}/2 \pi \mathbb{Z}$. Since the equations of the model are written in a conservative form, the natural framework to be used for the spatial discretization is the finite volume framework.  We hence introduce the control volume defined for $k \in \{1,...,N_\theta\}$
\begin{equation}
V_k = (\theta_{k-\frac{1}{2}},\theta_{k+\frac{1}{2}}).
\end{equation}
Let $\rho_k^n$ (resp. $u^{n}_{k+\frac{1}{2}}$) be the approximated value of the exact solution $\rho(t^n,\theta_k)$ (resp. $u(t^{n},\theta_{k+\frac{1}{2}})$), the classical upwind scheme for \eqref{1Dpolar} reads as
\begin{eqnarray}
\frac{\rho_k^{n+1} - \rho_k^n}{\Delta t} = \frac{\mathcal{F}_{k+\frac{1}{2}} - \mathcal{F}_{k-\frac{1}{2}}}{\Delta \theta}, \quad k \in \{1,...,N_\theta\},
\end{eqnarray}
where the numerical flux $\mathcal{F}_{k+\frac{1}{2}}$ and $\mathcal{F}_{k-\frac{1}{2}}$ are defined by
\begin{eqnarray*}
\mathcal{F}_{k+\frac{1}{2}} = \frac{\rho_{k+1}^{n+1} -  \rho_k^{n+1}}{\Delta \theta}- \textit{A}^{up}(u_{k+\frac{1}{2}}^{n+1},\rho_k^{n+1},\rho_{k+1}^{n+1}), \\
\mathcal{F}_{k-\frac{1}{2}} = \frac{\rho_{k}^{n+1} -  \rho_{k-1}^{n+1}}{\Delta \theta}- \textit{A}^{up}(u_{k-\frac{1}{2}}^{n+1},\rho_{k-1}^{n+1},\rho_k^{n+1}),
\end{eqnarray*}
with the advection numerical flux is given by
\begin{equation}\label{Auppolar}
\textit{A}^{up}(u,x_-,x_+) = \begin{cases}
u \, x_-, \quad \mbox{ if } u > 0, \\ 
u \, x_+, \quad \mbox{ if } u < 0.
\end{cases}
\end{equation}
The periodic flux condition on boundary reads as 
$$\mathcal{F}_{ \frac{1}{2}} = \mathcal{F}_{N_\theta + \frac{1}{2}} = \frac{\rho_{1}^{n+1} -  \rho_{N_\theta}^{n+1}}{\Delta \theta}- \textit{A}^{up}(u_{\frac{1}{2}}^{n+1},\rho_{N_\theta}^{n+1},\rho_1^{n+1})$$ and we recall that $u$ is periodic so we set the value $u^n_{ \frac{1}{2}} = u^n_{N_\theta + \frac{1}{2}}$.
The diffusion and advection terms are both treated implicitly, the scheme is then unconditionally stable.
We define the column vector $\rho^{n} = \begin{pmatrix} \rho_1^n & \rho_2^n & \dots & \rho_{N_\theta}^n  \end{pmatrix}^T$. As usual, see e.g. \cite{Allaire}, the discrete heat matrix $A \in M_{N_\theta} (\mathbb{R})$ with periodic flux condition on the boundary is defined as 
\begin{equation}\label{defApolar}
A = 
\begin{pmatrix} 
2 & -1 & & & -1 \\ -1 & 2 & \ddots \\ & \ddots & \ddots & \ddots \\& & \ddots & 2 & -1 \\ -1 & & & -1 & 2
\end{pmatrix}.
\end{equation}
Periodic flux condition adds the top right term and the bottom left term. Next, in order to use $A^{up}$ defined by equation \eqref{Auppolar}, we define $(u)^+=\max(u,0)$ and $(u)^-=\min(u,0)$ so $\textit{A}^{up}(u,\rho_k,\rho_{k+1}) = (u)^+ \rho_k + (u)^- \rho_{k+1}$. The discrete advection matrix $B^n \in M_{N_\theta} (\mathbb{R})$ with periodic flux condition on the boundary is then defined as in \cite{Allaire}
\begin{equation}\label{defBpolar}
B^n = 
\begin{pmatrix} 
\left(u_{\frac{3}{2}}^n\right)^+ & \left(u_{\frac{3}{2}}^n\right)^- & & &   \\ & \ddots & \ddots  \\ & & \left(u_{j+\frac{1}{2}}^n\right)^+ & \left(u_{j+\frac{1}{2}}^n\right)^-  \\ & & & \ddots & \left(u_{N_\theta-\frac{1}{2}}^n\right)^- \\ \left(u_{N_\theta+\frac{1}{2}}^n\right)^- & & & & \left(u_{N_\theta+\frac{1}{2}}^n\right)^+ 
\end{pmatrix}
\end{equation}
\begin{equation*}
-
\begin{pmatrix} 
\left(u_{\frac{1}{2}}^n\right)^- & & & & \left(u_{\frac{1}{2}}^n\right)^+  \\ \left(u_{\frac{3}{2}}^n\right)^+ & \ddots &   \\ & \left(u_{j-\frac{1}{2}}^n\right)^+ & \left(u_{j-\frac{1}{2}}^n\right)^- &  \\ & & \ddots & \ddots &  \\  & & & \left(u_{N_\theta-\frac{1}{2}}^n\right)^+ & \left(u_{N_\theta-\frac{1}{2}}^n\right)^-
\end{pmatrix}.
\end{equation*}
At each time step we have
\begin{equation*}
\frac{\rho^{n+1}-\rho^n}{\Delta t}= - \frac{1}{\Delta \theta^2} A \rho^{n+1}- \frac{1}{\Delta \theta} B^{n+1}\rho^{n+1}.
\end{equation*}
We use a standard numerical method to invert the matrix $A+\Delta \theta \, B^{n+1}+\frac{\Delta \theta^2}{\Delta t} I_{N_\theta}$. Finally, at each time step we resolve
\begin{equation*}
\rho^{n+1}= \left(A+\Delta \theta \, B^{n+1}+\frac{\Delta \theta^2}{\Delta t} I_{N_\theta}\right)^{-1} \, \frac{\Delta \theta^2}{\Delta t} \, \rho^{n}.
\end{equation*}

\subsection{Two dimensional case: polar case}
Let $\Omega \subset \mathbb{R}^2$ be the cytoplasm domain, the domain is obviously an annulus where the molecular markers cannot enter in the nucleus of the cell. Let us first recall the model  on the annulus $\Omega = B(0,R_{max}) \setminus B(0,R_{min})$ with $\vec{n}_\bx$ the unit normal vector to $\Omega$ at point $\bx \in \partial \Omega$ (we note $C(0,R)$ the circle of center $(0,0)$ and radius $R$):
\begin{eqnarray}
& \partial_t \rho (t,\bx) = \nabla . \left( \nabla \rho (t,\bx) -  \rho (t,\bx) \, \bu (t,\bx) \right) \label{rho1rec},  \mbox{ in } \Omega, \\
& (\nabla \rho (t,\bx) -  \rho (t,\bx) \, \bu (t,\bx)).\vec{n}_\bx =  0 \label{rho2rec},  \mbox{ on } C(0,R_{max}),  \\
& (\nabla \rho (t,\bx) -  \rho (t,\bx) \, \bu (t,\bx)).\vec{n}_\bx =  0 \label{rho3rec},  \mbox{ on } C(0,R_{min}).
\end{eqnarray}
Since the domain $\Omega$ is assumed to be an annulus, it is appropriate to introduce polar coordinates $r$ and $\theta$. Let $\bx=(r \cos(\theta), r \sin(\theta)) \in \Omega$, we have the following equations on $\frac{1}{r} \tilde{\rho}(t,r,\theta) = \rho(t,\bx)$ with $(r,\theta) \in [R_{min},R_{max}] \times \mathbb{R}/2\pi\mathbb{Z}$:
\begin{eqnarray}
\partial_t \tilde{\rho}(t,r,\theta) & = & \partial_r \left(r \partial_r \left(\frac{\tilde{\rho}(t,r,\theta)}{r}\right) - \tilde{\rho}(t,r,\theta) \bu_r (t,r,\theta) \right) \nonumber \\
& + &  \partial_{\theta} \left( \frac{1}{r^2} (\partial_{\theta}\tilde{\rho}(t,r,\theta) - \tilde{\rho}(t,r,\theta) \bu_\theta (t,r,\theta) ) \right) \label{rho1polar}, \mbox{ in } \Omega, \\
0 & = &  r \partial_r \left(\frac{\tilde{\rho}(t,r,\theta)}{r}\right) - \tilde{\rho}(t,r,\theta) \bu_r (t,r,\theta) \label{rho2polar}, \mbox{ on } C(0,R_{max}),  \\
0 &=& r \partial_r \left(\frac{\tilde{\rho}(t,r,\theta)}{r}\right) - \tilde{\rho}(t,r,\theta) \bu_r (t,r,\theta) \label{rho3polar}, \mbox{ on } C(0,R_{min}).
\end{eqnarray}
We recall the potential case $\bu(t,\bx) = \nabla c(t,\bx)$. Laplace equation on $c$ with non appropriate Neumann conditions on a bounded domain is ill-posed,  see \cite{Allaire} e.g. In order to handle this problem, we add a degradation term:
\begin{eqnarray}
& -\Delta c(\bx) + \alpha \, c(\bx) = 0 \label{c1rec}, \mbox{ in } \Omega,\\
& \nabla c(\bx).\vec{n}_\bx = \rho(t,\bx) \label{c2rec}, \mbox{ on } C(0,R_{max}), \\
& \nabla c(\bx).\vec{n}_\bx = 0 \label{c3rec}, \mbox{ on } C(0,R_{min}).
\end{eqnarray}
We also give the previous equations in polar coordinates by $\frac{1}{r} \tilde{c}(r,\theta) = c(\bx)$ with $(r,\theta) \in [R_{min},R_{max}] \times \mathbb{R}/2\pi\mathbb{Z}$:
\begin{eqnarray}
& - \partial_r \left(r \partial_r \left(\frac{\tilde{c}(r,\theta)}{r}\right) \right)- \frac{1}{r^2} \partial_{\theta \theta} \tilde{c} (r,\theta) + \alpha \, \tilde{c} (r,\theta) = 0, \mbox{ in } \Omega \label{c1polar},\\
& \partial_r \left(\frac{\tilde{c} (r,\theta)}{r}\right) =  \frac{\tilde{\rho} (R_{max},\theta)}{R_{max}}, \mbox{ on } C(0,R_{max}) \label{c2polar}, \\
& \partial_r \left(\frac{\tilde{c} (r,\theta)}{r}\right) = 0, \mbox{ on } C(0,R_{min}) \label{c3polar}.
\end{eqnarray}
If we consider dynamical exchange of markers at the active boundary, for $\bx \in \Gamma=C(0,R_{max})$ we have the evolution in time of $\mu(t,\bx)$.
\begin{eqnarray}\label{mupolar}
\partial_t \mu(t,\bx) =\rho(t,\bx) -\mu(t,\bx) , \mbox{ on } C(0,R_{max}).
\end{eqnarray}
We replace then \eqref{rho2polar} by
\begin{eqnarray}
-\partial_t \mu(t,\theta) & = &  r \partial_r \left(\frac{\tilde{\rho}(t,r,\theta)}{r}\right) - \tilde{\rho}(t,r,\theta) \bu_r (t,r,\theta) , \mbox{ on } C(0,R_{max}),
\end{eqnarray}
and \eqref{c2polar} by
\begin{eqnarray}
\partial_r \left(\frac{\tilde{c} (r,\theta)}{r}\right) =  \, \mu(t,\theta), \mbox{ on } C(0,R_{max}).
\end{eqnarray}
The transversal case is then $$\bu(t,r,\theta) = \frac{\tilde{\rho}(t,R_{max},\theta)}{R_{max}} \be_r,$$
and we can also write the transversal case for the dynamical exchange model
$$\bu(t,r,\theta) = \mu(t,\theta) \be_r,$$

Let $t^n = n\, \Delta t$ be the time discretization and $\{r_j=R_{min} + j \, \Delta r, j \in \{1,...,N_r\}\}$ (resp. $\{\theta_k=k \, \Delta \theta, k \in \{1,...,N_\theta\}\}$) be the space discretization of the bounded interval $[R_{min},R_{max}]$ (resp. periodic interval $\mathbb{R}/ 2 \pi \mathbb{Z}$). We introduce the control volume $W_{(j,k)}  \subset \mathbb{R}^2$
\begin{equation*}
W_{(j,k)} = \left(r_{j-\frac{1}{2}},r_{j+\frac{1}{2}}\right) \times \left(\theta_{k-\frac{1}{2}},\theta_{k+\frac{1}{2}}\right).
\end{equation*}

Let $\tilde{P}^{n}_{(j,k)}$ (resp. $\mu^n_{k}$) be the approximated value of the exact solution $\tilde{\rho}(t^n,r_j,\theta_k)$ (resp. $\mu(t^n,\theta_k)$). 
Let  $\tilde{c}_{(j,k)}$ be the approximated value of the exact solution $\tilde{c}(r_j,\theta_k)$. 

\subsubsection{Equation on $\mu$}
In the dynamical exchange model, we can resolve at each time step the discretization of equation \eqref{mupolar} for $k \in \{1,...,N_y\}$
\begin{eqnarray*}
\mu^{n+1}_k=\mu^{n}_k + \Delta t \, (\rho^n_k-\mu^n_k).
\end{eqnarray*}

\subsubsection{Equation on $\tilde{c}$}
For simplicity, we call $\mathcal{F}$ the numerical flux as in the 1D case, we can write the following scheme for equation \eqref{c1polar} for $(j,k) \in \{1,...,N_r\} \times \{1,...,N_\theta\}$
\begin{eqnarray*}
- \left(\frac{\mathcal{F}_{(j+\frac{1}{2},k)} - \mathcal{F}_{(j-\frac{1}{2},k)}}{\Delta r} +  \frac{\mathcal{F}_{(j,k+\frac{1}{2})} - \mathcal{F}_{(j,k-\frac{1}{2})}}{\Delta \theta}\right) + \alpha \, \tilde{c}_{(j,k)} = 0.
\end{eqnarray*}
In order to use finite volume we define the numerical flux by
\begin{eqnarray*}
\mathcal{F}_{(j+\frac{1}{2},k)} = r_{j+\frac{1}{2}}\, \frac{\frac{\tilde{c}_{(j+1,k)}}{r_{j+1}} - \frac{\tilde{c}_{(j,k)}}{r_{j}} }{\Delta r}, & \quad &
\mathcal{F}_{(j-\frac{1}{2},k)} = r_{j-\frac{1}{2}}\, \frac{\frac{\tilde{c}_{(j,k)}}{r_{j}} - \frac{\tilde{c}_{(j-1,k)}}{r_{j-1}} }{\Delta r},\\
\mathcal{F}_{(j,k+\frac{1}{2})} = \frac{1}{r_{j}^2} \frac{\tilde{c}_{(j,k+1)} - \tilde{c}_{(j,k)} }{\Delta \theta}, & \quad &
\mathcal{F}_{(j,k-\frac{1}{2})} = \frac{1}{r_{j}^2} \frac{\tilde{c}_{(j,k)} - \tilde{c}_{(j,k-1)} }{\Delta \theta}.
\end{eqnarray*}
The zero flux boundary condition \eqref{c3polar} impose that $\mathcal{F}_{(\frac{1}{2},k)} = 0$ and the boundary condition \eqref{c2polar} $\mathcal{F}_{(N_r + \frac{1}{2},k)} = r_{N_r+\frac{1}{2}} \, \mu^n_k$ for $k \in \{1,...,N_\theta\}$. Similarly, the periodic conditions impose for $j \in \{1,...,N_r\}$ 
$$\mathcal{F}_{(j,N_\theta+\frac{1}{2})} = \mathcal{F}_{(j,\frac{1}{2})}=\frac{1}{r_{j}^2} \frac{\tilde{c}_{(j,1)} - \tilde{c}_{(j,N_\theta)} }{\Delta \theta}.$$
We define the column vector $\mathcal{C}$ by $\mathcal{C}(k+(j-1) N_\theta) = \tilde{c}_{(j,k)}$ with $(j,k) \in \{1,...,N_r\} \times \{1,...,N_\theta\}$:
\begin{equation*}
\mathcal{C} =
\left(
\tilde{c}_{(1,1)} \, \dots \,  \tilde{c}_{(1,N_\theta)} \, \tilde{c}_{(2,1)} \, \dots \, \tilde{c}_{(2,N_\theta)} \, \dots  \, \tilde{c}_{(N_r,N_\theta)}
\right)^T
\end{equation*}
For $\Delta r=\Delta \theta$ the rigidity matrix $\mathcal{A}$ is defined by
\begin{eqnarray}\label{A2D}
\mathcal{A} = 
\begin{pmatrix} 

 \\
& \ddots & \ddots & \ddots &  \\ 
& & -\frac{r_{j-\frac{1}{2}}}{r_{j-1}} \,Id & \frac{r_{j-\frac{1}{2}} + r_{j+\frac{1}{2}}}{r_j} \, Id & -\frac{r_{j+\frac{1}{2}}}{r_{j+1}} \, Id \\ & & & \ddots & \ddots & \ddots &  \\

\end{pmatrix} \nonumber \\
+ 
\begin{pmatrix} 
\frac{1}{r_{1}^2} A  \\  & \frac{1}{r_{2}^2} A  \\  & & \ddots \\ & & & \frac{1}{r_{N_r-1}^2} A \\ & & & & \frac{1}{r_{N_r}^2} A
\end{pmatrix}.
\end{eqnarray}
The flux boundary condition $C(0,R_{max})$ imposes this right hand side column vector of length $N_r \, N_\theta$:
\begin{equation*}
\mathcal{R}^n_c = r_{N_r+\frac{1}{2}} \begin{pmatrix} 0 \\ \vdots \\ 0 \\ \left(\frac{\tilde{\rho}^n_{N_r,k}}{r_{N_r}}\right)_{k} \end{pmatrix} \mbox{ or in the exchange case } \mathcal{R}^n_c = r_{N_r+\frac{1}{2}} \begin{pmatrix} 0 \\ \vdots \\ 0 \\ (\mu^n_k)_{k} \end{pmatrix}.
\end{equation*}
We use a standard numerical method to invert the symmetric positive definite matrix $\frac{1}{\Delta r^2} \mathcal{A}+ \alpha I_{N_r N_\theta} $ and then resolve at each time step
\begin{equation*}
\mathcal{C} = \left(\frac{1}{\Delta r^2}\mathcal{A}
+ \alpha I_{N_r N_\theta} \right)^{-1} \, \frac{1}{\Delta r} \mathcal{R}^n_c.
\end{equation*}

\subsubsection{Equation on $\tilde{\rho}$}
For simplicity, we call $\mathcal{F}$ the numerical flux as in the 1D case, we can write the following scheme for equation \eqref{rho1polar}: for $(j,k) \in \{1,...,N_r\} \times \{1,...,N_\theta\}$
\begin{eqnarray*}
\frac{\tilde{P}_{(j,k)}^{n+1} - \tilde{P}_{(j,k)}^n}{\Delta t} = \frac{\mathcal{F}_{(j+\frac{1}{2},k)} - \mathcal{F}_{(j-\frac{1}{2},k)}}{\Delta r} +  \frac{\mathcal{F}_{(j,k+\frac{1}{2})} - \mathcal{F}_{(j,k-\frac{1}{2})}}{\Delta \theta}.
\end{eqnarray*}
We define the numerical flux by
\begin{eqnarray*}
\mathcal{F}_{(j+\frac{1}{2},k)} = r_{j+\frac{1}{2}}\, \frac{\frac{\tilde{P}^{n+1}_{(j+1,k)}}{r_{j+1}} - \frac{\tilde{P}^{n+1}_{(j,k)}}{r_{j}} }{\Delta r}  - A^{up} \left(u^{n+1}_{(j+\frac{1}{2},k)},\tilde{P}^{n+1}_{(j,k)},\tilde{P}^{n+1}_{(j+1,k)}\right),\\
\mathcal{F}_{(j-\frac{1}{2},k)} = r_{j-\frac{1}{2}}\, \frac{\frac{\tilde{P}^{n+1}_{(j,k)}}{r_{j}} - \frac{\tilde{P}^{n+1}_{(j-1,k)}}{r_{j-1}} }{\Delta r}  - A^{up} \left(u^{n+1}_{(j-\frac{1}{2},k)},\tilde{P}^{n+1}_{(j-1,k)},\tilde{P}^{n+1}_{(j,k)}\right),\\
\mathcal{F}_{(j,k+\frac{1}{2})}  = \frac{1}{r_{j}^2} \, \left( \frac{\tilde{P}^{n+1}_{(j,k+1)} - \tilde{P}^{n+1}_{(j,k)} }{\Delta \theta} - A^{up} \left(u^{n+1}_{(j,k+\frac{1}{2})},\tilde{P}^{n+1}_{(j,k)},\tilde{P}^{n+1}_{(j,k+1)}\right) \right),\\
\mathcal{F}_{(j,k+\frac{1}{2})}=  \frac{1}{r_{j}^2} \, \left(\frac{\tilde{P}^{n+1}_{(j,k)} - \tilde{P}^{n+1}_{(j,k-1)} }{\Delta \theta} - A^{up} \left(u^{n+1}_{(j,k-\frac{1}{2})},\tilde{P}^{n+1}_{(j,k-1)},\tilde{P}^{n+1}_{(j,k)}\right)\right).
\end{eqnarray*}
In the transversal case \eqref{MT}, we take $\bu_r = \frac{\tilde{\rho}(R_{max},\theta)}{r}$ and $\bu_\theta = 0$ we define at time $t^n$
\begin{eqnarray*}
u^n_{(j+\frac{1}{2},k)}= - \frac{\tilde{\rho}^{n}_{(N_r,k)}}{r_{N_r}}  \mbox{ or in} & \mbox{ the } & \mbox{exchange case } u^n_{(j+\frac{1}{2},k)}= - \mu_k^n, \\
u^n_{(j-\frac{1}{2},k)}= - \frac{\tilde{\rho}^{n}_{(N_r,k)}}{r_{N_r}} \mbox{ or in} & \mbox{ the } & \mbox{exchange case } u^n_{(j-\frac{1}{2},k)} = - \mu_k^n, \\
u^n_{(j,k+\frac{1}{2})}= 0, & \quad &
u^n_{(j,k-\frac{1}{2})}= 0.
\end{eqnarray*}
In the potential case \eqref{c}, we have $\bu_r = \partial_r \left(\frac{\tilde{c}}{r}\right)$ and $\bu_\theta = \frac{1}{r}\partial_\theta \tilde{c}$ we define at time $t^n$
\begin{eqnarray*}
u^n_{(j+\frac{1}{2},k)}=\frac{\frac{\tilde{c}_{(j+1,k)}}{r_{j+1}} - \frac{\tilde{c}_{(j,k)}}{r_{j}} }{\Delta r}, & \quad &
u^n_{(j-\frac{1}{2},k)}=\frac{\frac{\tilde{c}_{(j,k)}}{r_{j}} - \frac{\tilde{c}_{(j-1,k)}}{r_{j-1}}}{\Delta r}, \\
u^n_{(j,k+\frac{1}{2})}=\frac{1}{r_j}\frac{\tilde{c}_{(j,k+1)} - \tilde{c}_{(j,k)} }{\Delta \theta}, & \quad &
u^n_{(j,k-\frac{1}{2})}=\frac{1}{r_j}\frac{\tilde{c}_{(j,k)} - \tilde{c}_{(j,k-1)} }{\Delta \theta}.
\end{eqnarray*}
The zero flux boundary condition \eqref{rho3polar} impose that $\mathcal{F}_{(\frac{1}{2},k)} = 0$. In the simplified model, the boundary condition \eqref{rho2polar} $\mathcal{F}_{(N_r + \frac{1}{2},k)} = 0$ for $k \in \{1,...,N_\theta\}$ and in the model with exchange, we have $\mathcal{F}_{(N_r + \frac{1}{2},k)} = - \frac{\mu^{n+1}_k-\mu^n_k}{\Delta t}$ for $k \in \{1,...,N_\theta\}$. Similarly, the periodic conditions impose for $j \in \{1,...,N_r\}$ 
$$\mathcal{F}_{(j,N_\theta+\frac{1}{2})} = \mathcal{F}_{(j,\frac{1}{2})}=\frac{1}{r_{j}^2} \, \left(\frac{\tilde{P}^{n+1}_{(j,1)} - \tilde{P}^{n+1}_{(j,N_\theta)} }{\Delta \theta} - A^{up} \left(u^{n+1}_{(j,\frac{1}{2})},\tilde{P}^{n+1}_{(j,N_\theta)},\tilde{P}^{n+1}_{(j,1)}\right)\right).$$
We define the column vector $\mathcal{P}^n$ by $\mathcal{P}^n(k+(j-1) N_\theta) = \tilde{P}^n_{(j,k)}$ with $(j,k) \in \{1,...,N_r\} \times \{1,...,N_\theta\}$:
\begin{equation*}
\mathcal{P}^n =
\left( 
\tilde{P}^n_{(1,1)} \, \dots  \, \tilde{P}^n_{(1,N_\theta)} \,  \tilde{P}^n_{(2,1)} \, \dots \, \tilde{P}^n_{(2,N_\theta)} \, \dots  \, \tilde{P}^n_{(N_r,N_\theta)}
\right)^T.
\end{equation*}
We define the following diagonal matrices for $j \in \{1,...,N_r\}$,  $U^{+}_{j+\frac{1}{2}} \in M_{N_\theta}(\mathbb{R})$ and $U^{-}_{j+\frac{1}{2}} \in M_{N_\theta}(\mathbb{R})$:
\begin{eqnarray*}
U^{+}_{j+\frac{1}{2}} = 
\begin{pmatrix} 
\ddots & & & &  \\ & (u_{(j+\frac{1}{2},k-1)}^n)^+ &   \\ & & (u_{(j+\frac{1}{2},k)}^n)^+  &  &  \\ & & & (u_{(j+\frac{1}{2},k+1)}^n)^+ &   \\  & & & & \ddots
\end{pmatrix}, \\
U^{-}_{j+\frac{1}{2}} = 
\begin{pmatrix} 
\ddots & & & &  \\ & (u_{(j+\frac{1}{2},k-1)}^n)^- &   \\ & & (u_{(j+\frac{1}{2},k)}^n)^-  &  &  \\ & & & (u_{(j+\frac{1}{2},k+1)}^n)^- &   \\  & & & & \ddots
\end{pmatrix}.
\end{eqnarray*}
Thus we can define:
\begin{eqnarray}\label{B2D}
\mathcal{B}^n & =& 
\begin{pmatrix} 
U_{\frac{3}{2}}^+ & U_{\frac{3}{2}}^- & & &   \\ & \ddots & \ddots  \\ & & U_{j+\frac{1}{2}}^+ & U_{j+\frac{1}{2}}^-  \\ & & & \ddots & U_{N_r-\frac{1}{2}}^- \\ & & & & 0
\end{pmatrix}
-
\begin{pmatrix} 
0 & & & &  \\ U_{\frac{3}{2}}^+ & \ddots &   \\ & U_{j-\frac{1}{2}}^+ & U_{j-\frac{1}{2}}^- &  \\ & & \ddots & \ddots &  \\  & & & U_{N_r-\frac{1}{2}}^+ & U_{N_r-\frac{1}{2}}^-
\end{pmatrix} \nonumber \\
& + &
\begin{pmatrix} 
\frac{1}{r_1^2} B^n  \\  & \frac{1}{r_2^2} B^n  \\  & & \ddots \\ & & & \frac{1}{r_{N_r-1}^2}  B^n \\ & & & &  \frac{1}{r_{N_r}^2} B^n
\end{pmatrix}
\end{eqnarray}
In the simplified model, we have at each time step
\begin{equation*}
\frac{\mathcal{P}^{n+1}-\mathcal{P}^{n}}{\Delta t} = -\frac{1}{\Delta r^2} \mathcal{A}\mathcal{P}^{n+1}- \frac{1}{\Delta r} \mathcal{B}^{n+1} \mathcal{P}^{n+1}.
\end{equation*}
In the exchange model, the flux boundary condition on $C(0,R_{max})$ imposes this right hand side column vector of length $N_r \, N_\theta$:
\begin{equation*}
\mathcal{R}_{exchange}^n = - \begin{pmatrix} 0 \\ \vdots \\ 0 \\ (\frac{\mu^{n+1}_k-\mu^n_k}{\Delta t})_{k} \end{pmatrix}.
\end{equation*}
We have at each time step
\begin{equation*}
\frac{\mathcal{P}^{n+1}-\mathcal{P}^{n}}{\Delta t} = -\frac{1}{\Delta r^2} \mathcal{A}\mathcal{P}^{n+1}- \frac{1}{\Delta r} \mathcal{B}^{n+1}  \mathcal{P}^{n+1} + \frac{1}{\Delta r} \mathcal{R}_{exchange}^n.
\end{equation*}
We use a standard numerical method to invert the matrix $\mathcal{A} + \Delta r \, \mathcal{B}^{n+1}+ \frac{\Delta r^2}{\Delta t} \, I_{N_r N_\theta}$ and then resolve at each time step
\begin{equation*}
\mathcal{P}^{n+1} = \left(\mathcal{A}
 + \Delta r \, \mathcal{B}^{n+1} + \frac{\Delta r^2}{\Delta t} \, I_{N_r N_\theta}\right)^{-1} \, \left(\frac{\Delta r^2}{\Delta t} \mathcal{P}^{n}+\mathcal{R}_{exchange}^n\right).
\end{equation*}

\subsection{Graphics}\label{Graphics}
\begin{center}
\begin{figure}[H]
\begin{tabular}{ll}
\includegraphics[scale=.25]{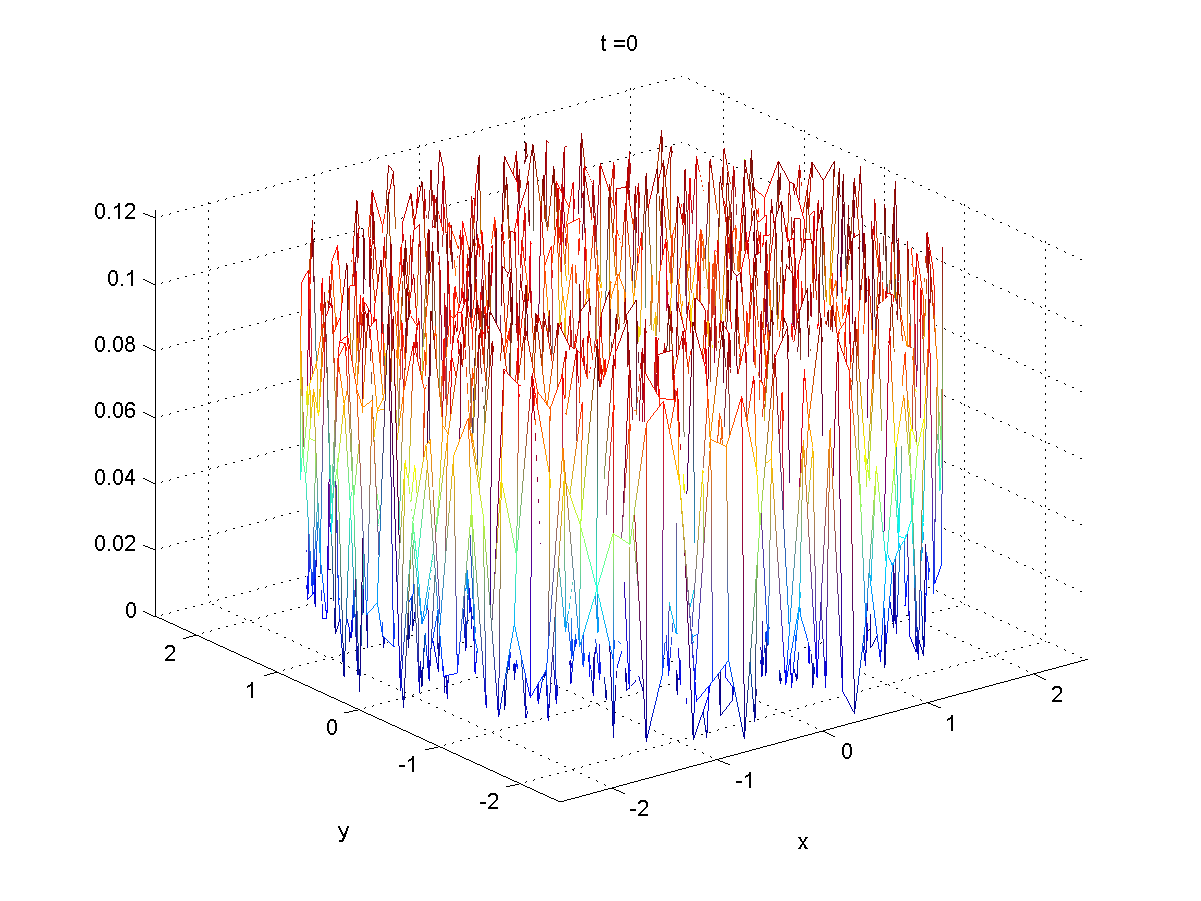} & \includegraphics[scale=.25]{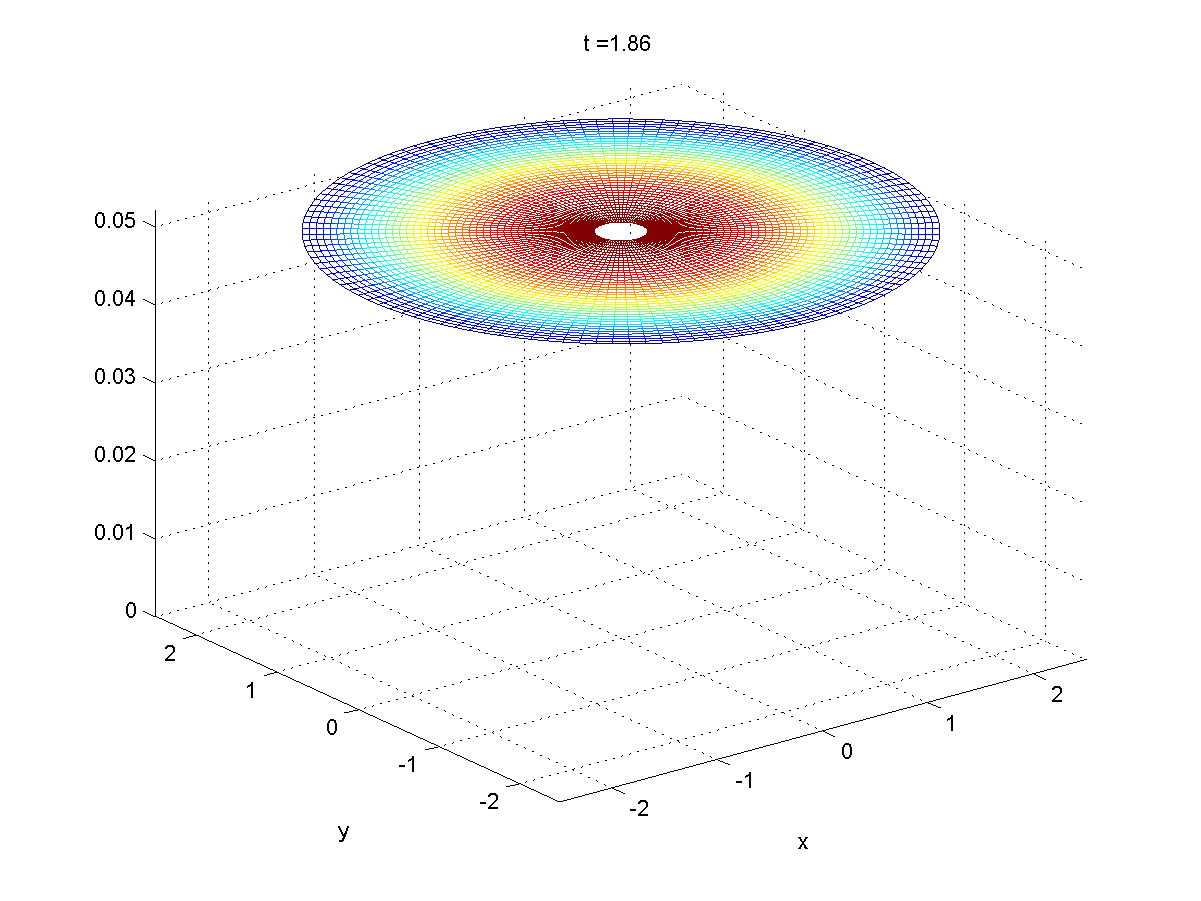}
\end{tabular}
 \caption{Numerical Simulations on $\Omega = [0.2,2.5] \times \mathbb{R}/2\pi \mathbb{Z}$ and all parameters equal to 1 with random initial conditions. Fo. For $M=0.01$ small, steady state is isotropic.}\label{figpolar2}
\end{figure}
\end{center}

We use the numerical analysis done in this article to implement it with Matlab. Simulations have first been done in the transversal case $\bu_T$. The following behaviour was obtained:
\begin{center}
\begin{figure}[H]
\begin{tabular}{ll}
 \includegraphics[scale=.25]{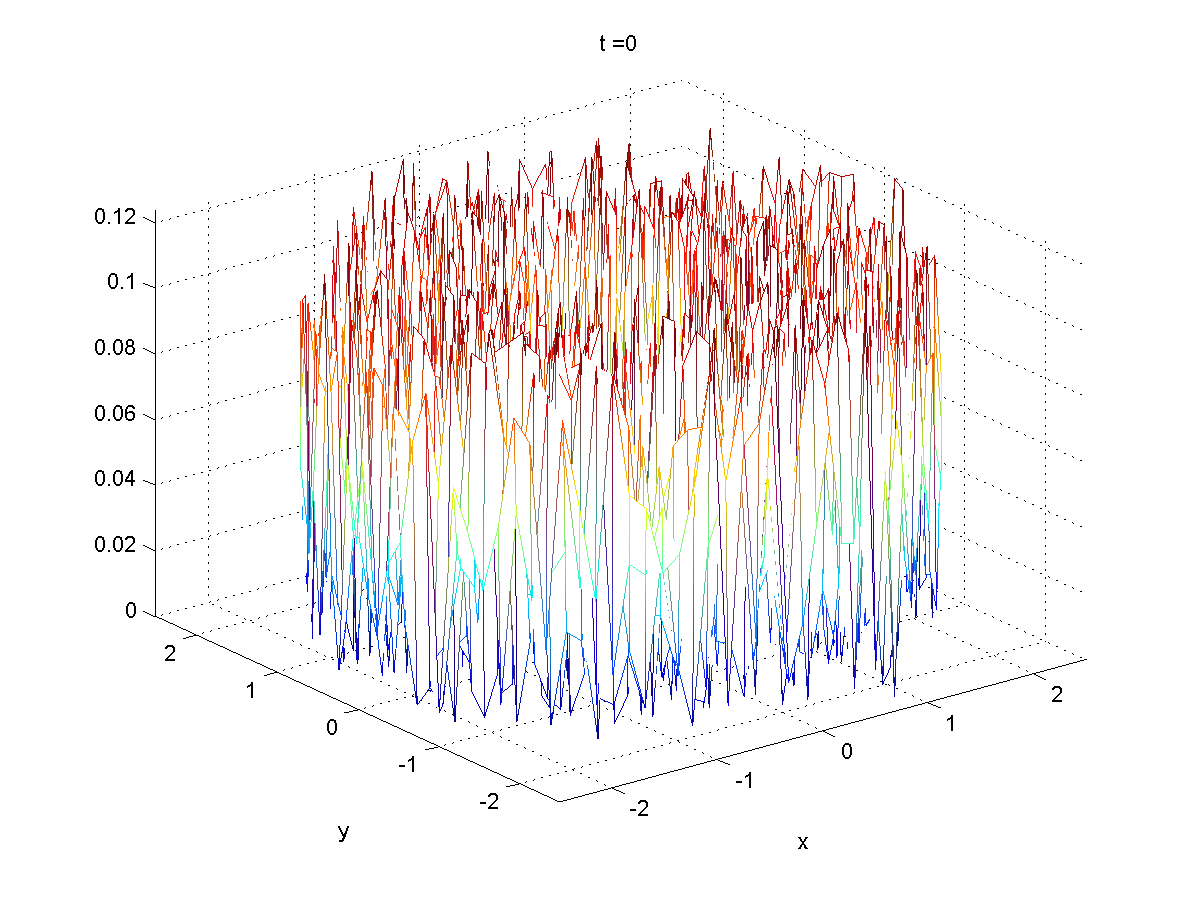}   & \includegraphics[scale=.25]{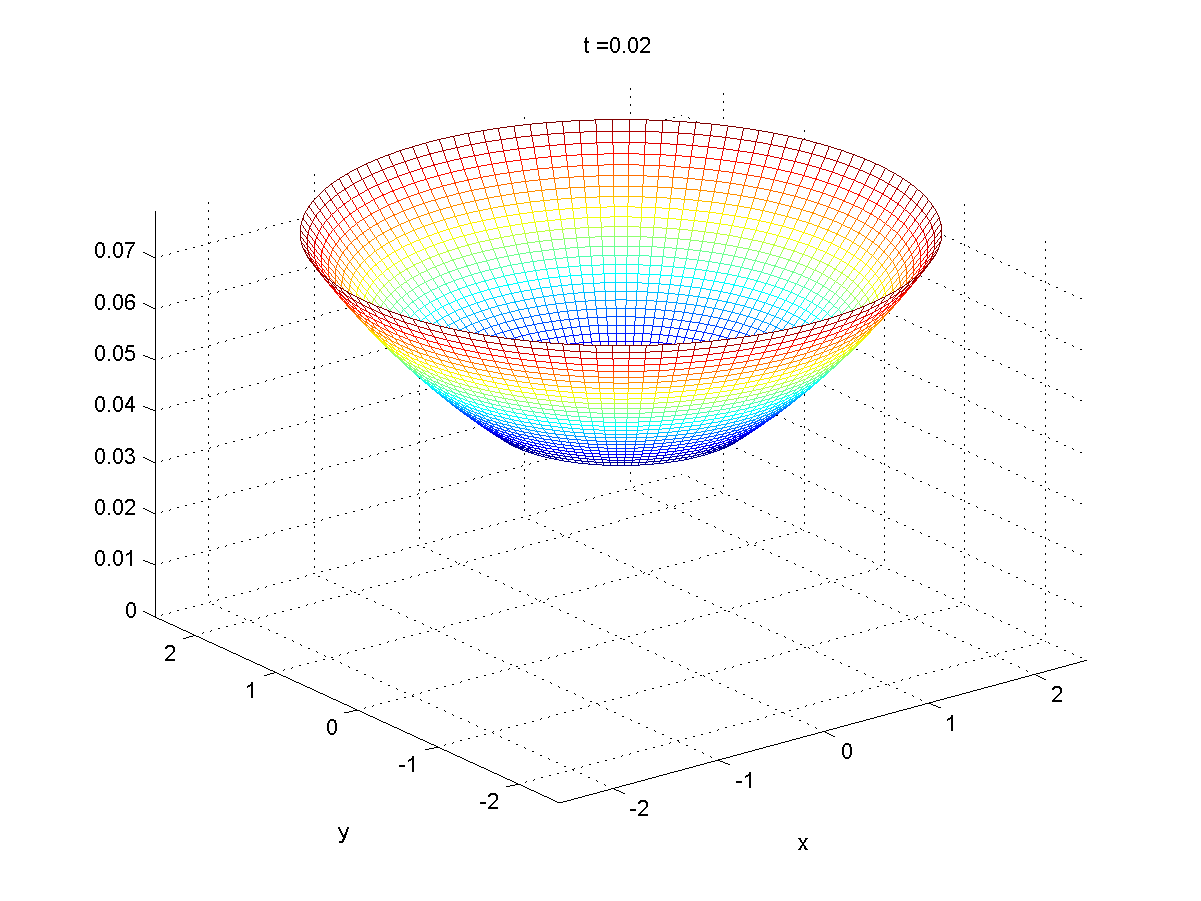} \\
  \includegraphics[scale=.25]{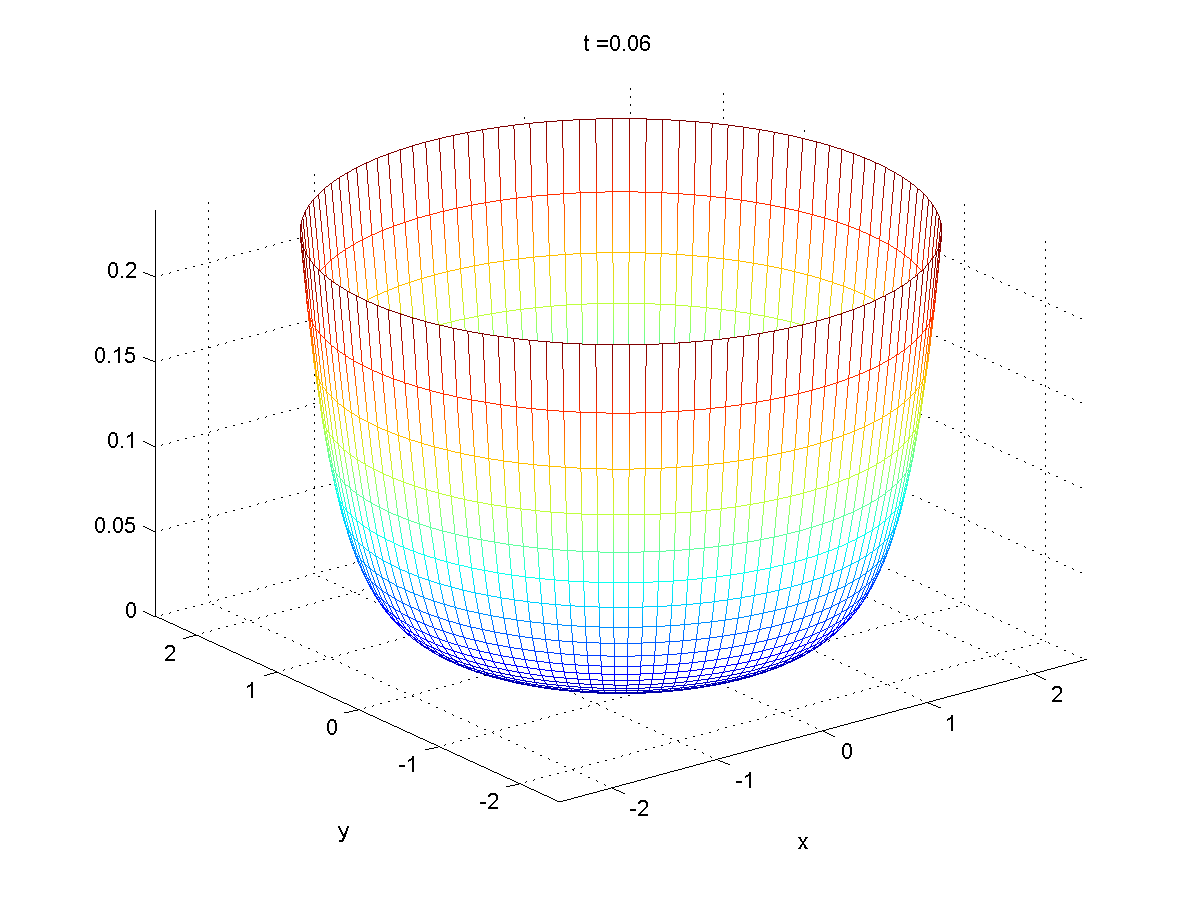}   & \includegraphics[scale=.25]{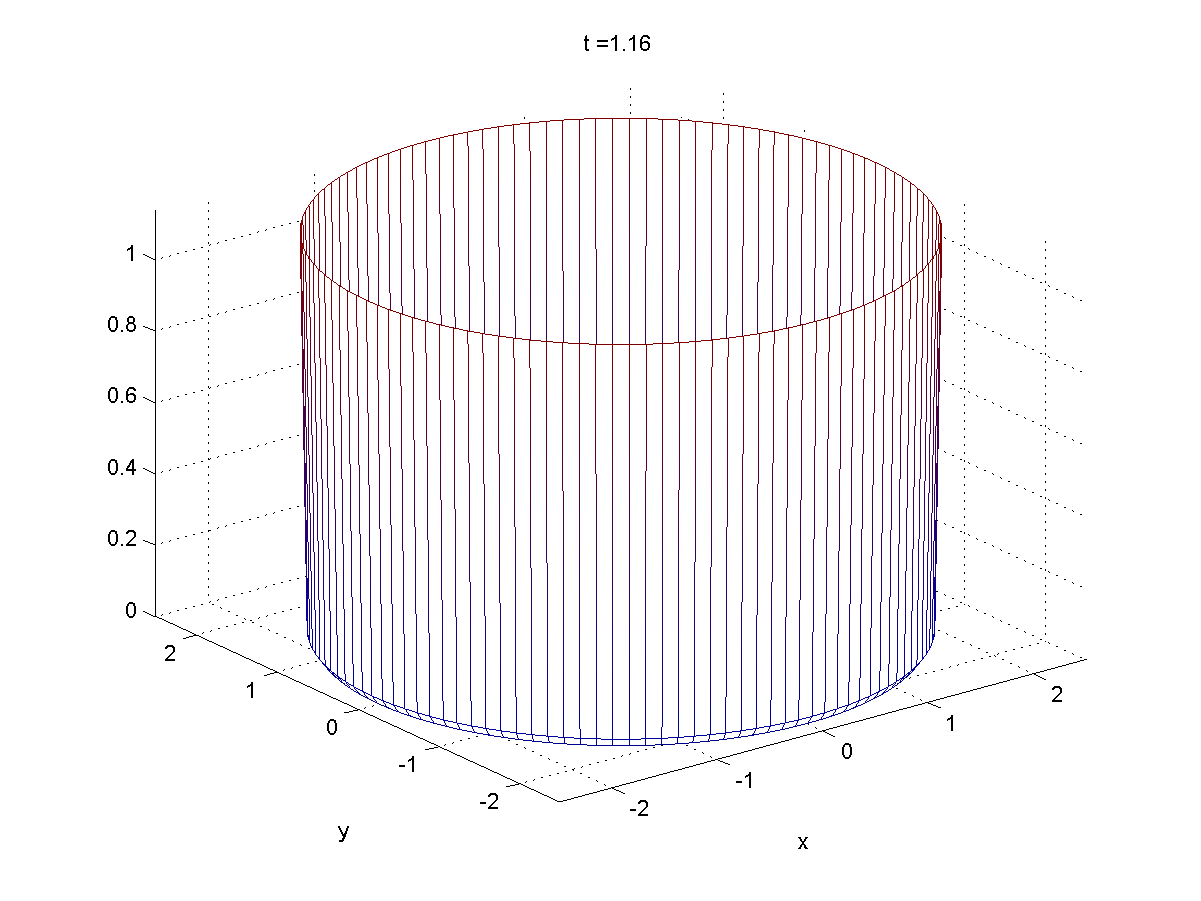}
\end{tabular} 
\caption{Numerical Simulations on $\Omega = [0.2,2.5] \times \mathbb{R}/2\pi \mathbb{Z}$ and all parameters equal to 1 with random initial conditions. For $M=10$ greater enough, no symmetry breaking appears, molecular markers are uniformly distributed on the membrane.}\label{figMT}
\end{figure}
\end{center} 
The dichotomy planned on $M$ by the heuristic holds true. Indeed we have done simulations for small $M$.
\begin{center}
\begin{figure}[H]
\begin{tabular}{ll}
\includegraphics[scale=.25]{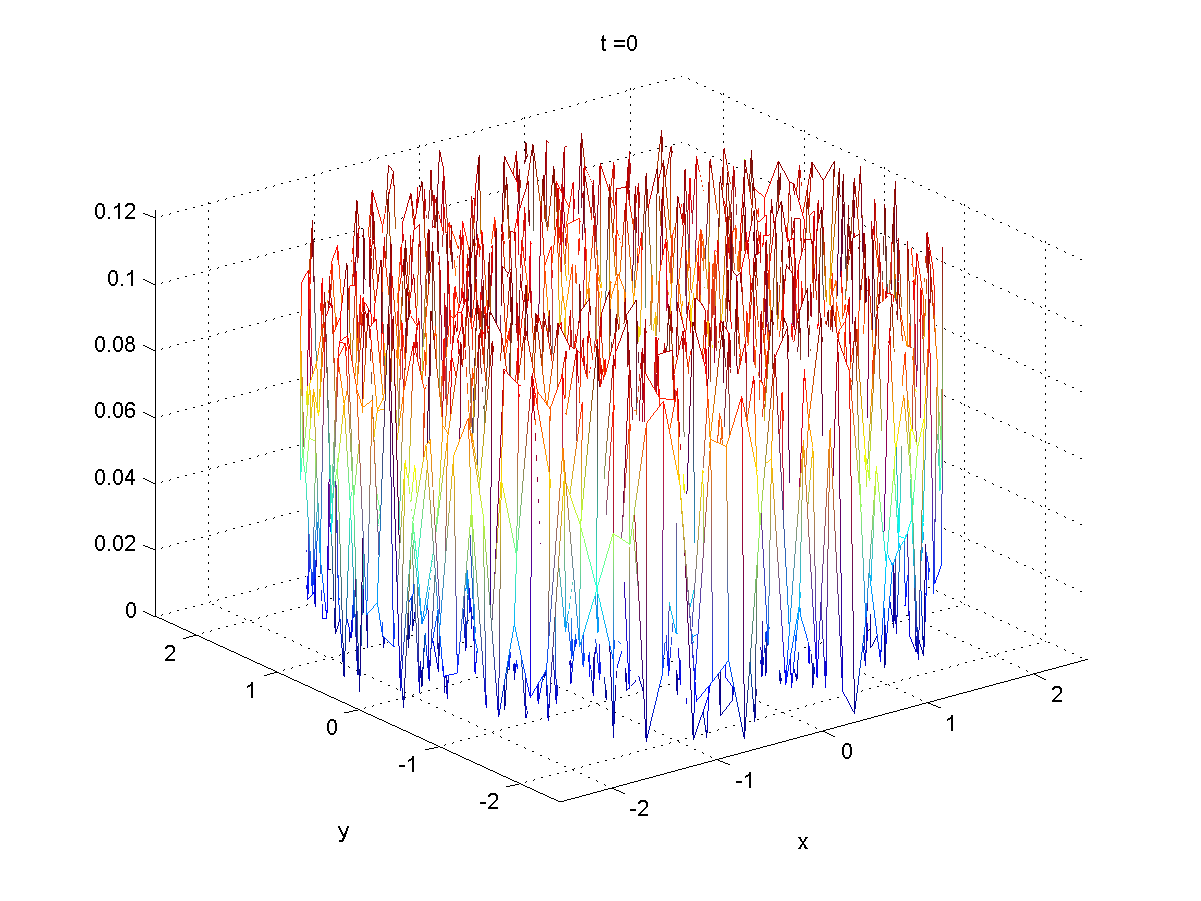} & \includegraphics[scale=.25]{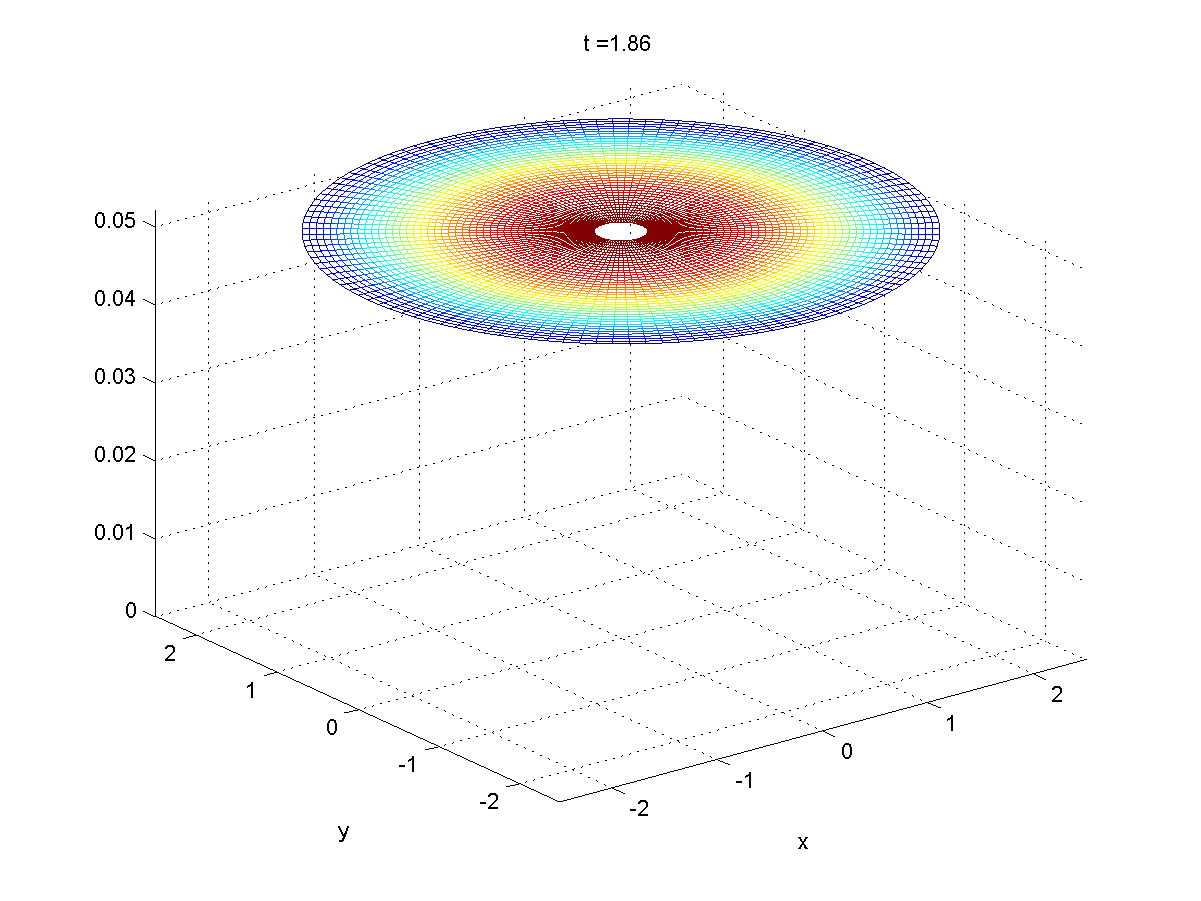}
\end{tabular}
 \caption{Numerical Simulations on $\Omega = [0.2,2.5] \times \mathbb{R}/2\pi \mathbb{Z}$ and all parameters equal to 1 with random initial conditions. Fo. For $M=0.01$ small, steady state is isotropic.}\label{figMT2}
\end{figure}
\end{center}
Then, simulations have been done in the potential case $\bu_P$:
\begin{center}
\begin{figure}[H]
\begin{tabular}{ll}
 \includegraphics[scale=.25]{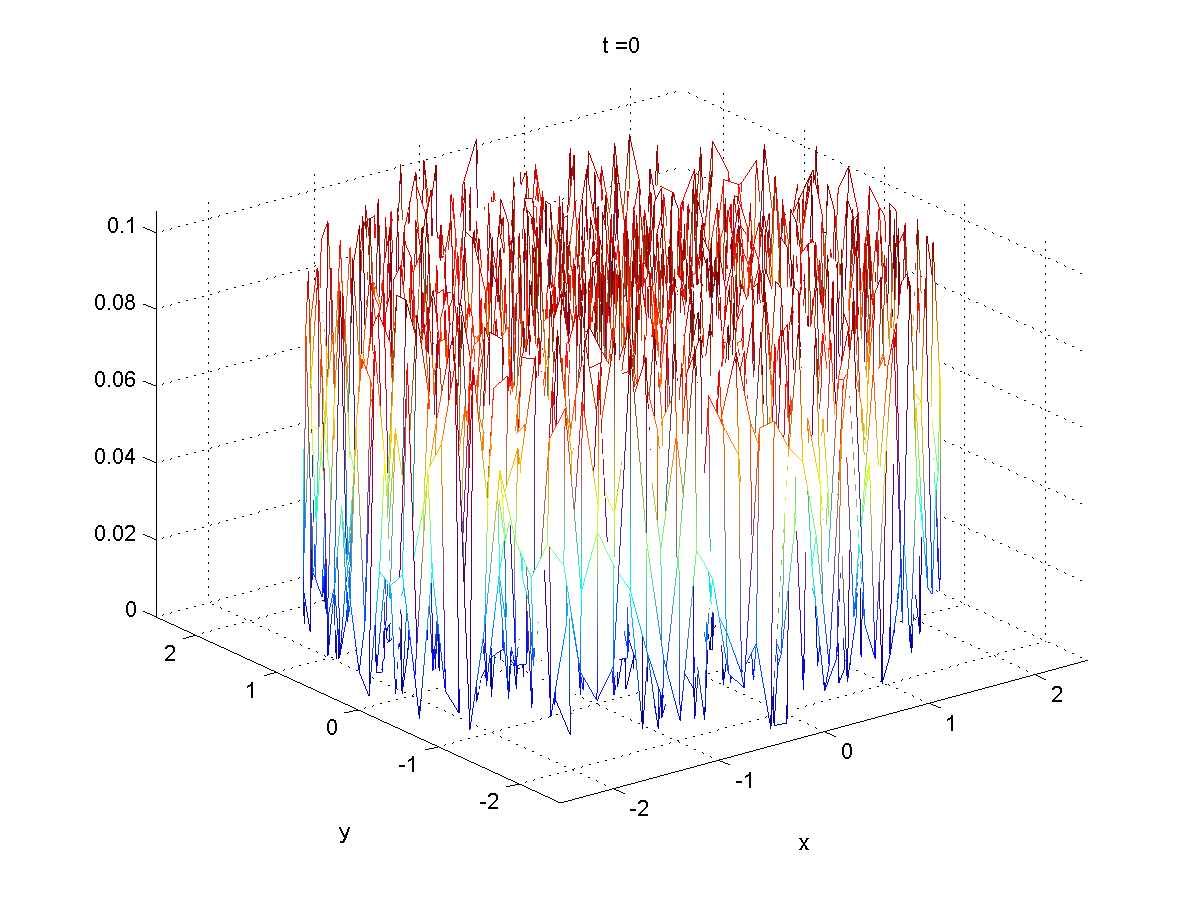}   & \includegraphics[scale=.25]{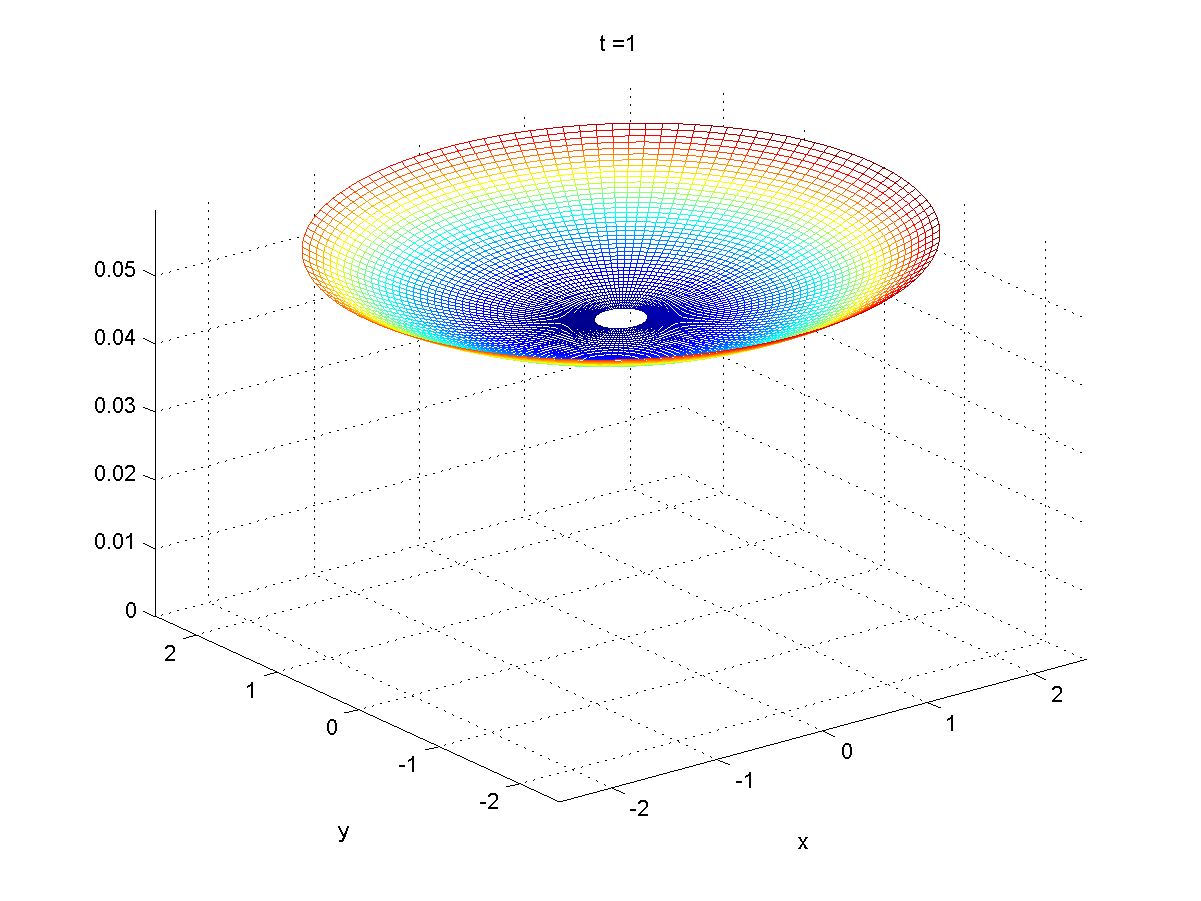} \\
  \includegraphics[scale=.25]{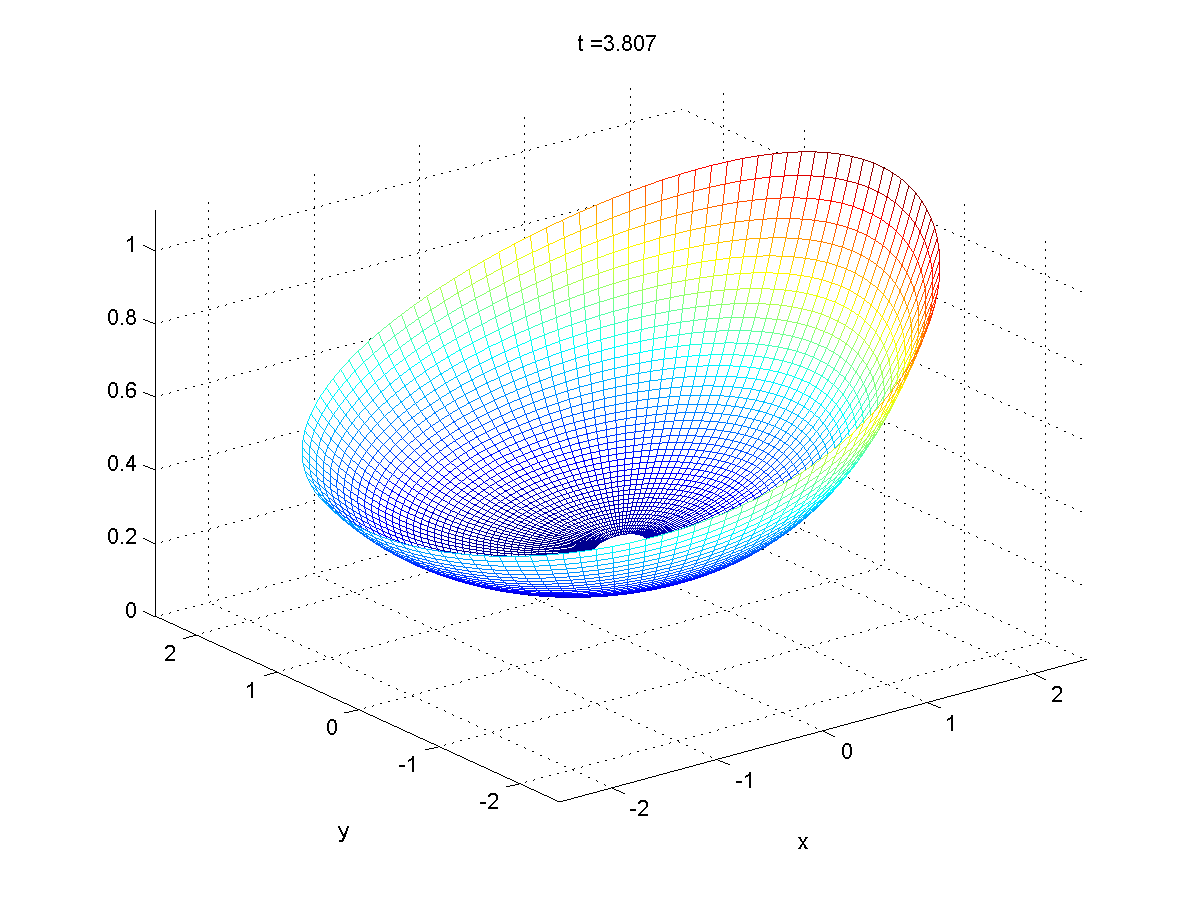}   & \includegraphics[scale=.25]{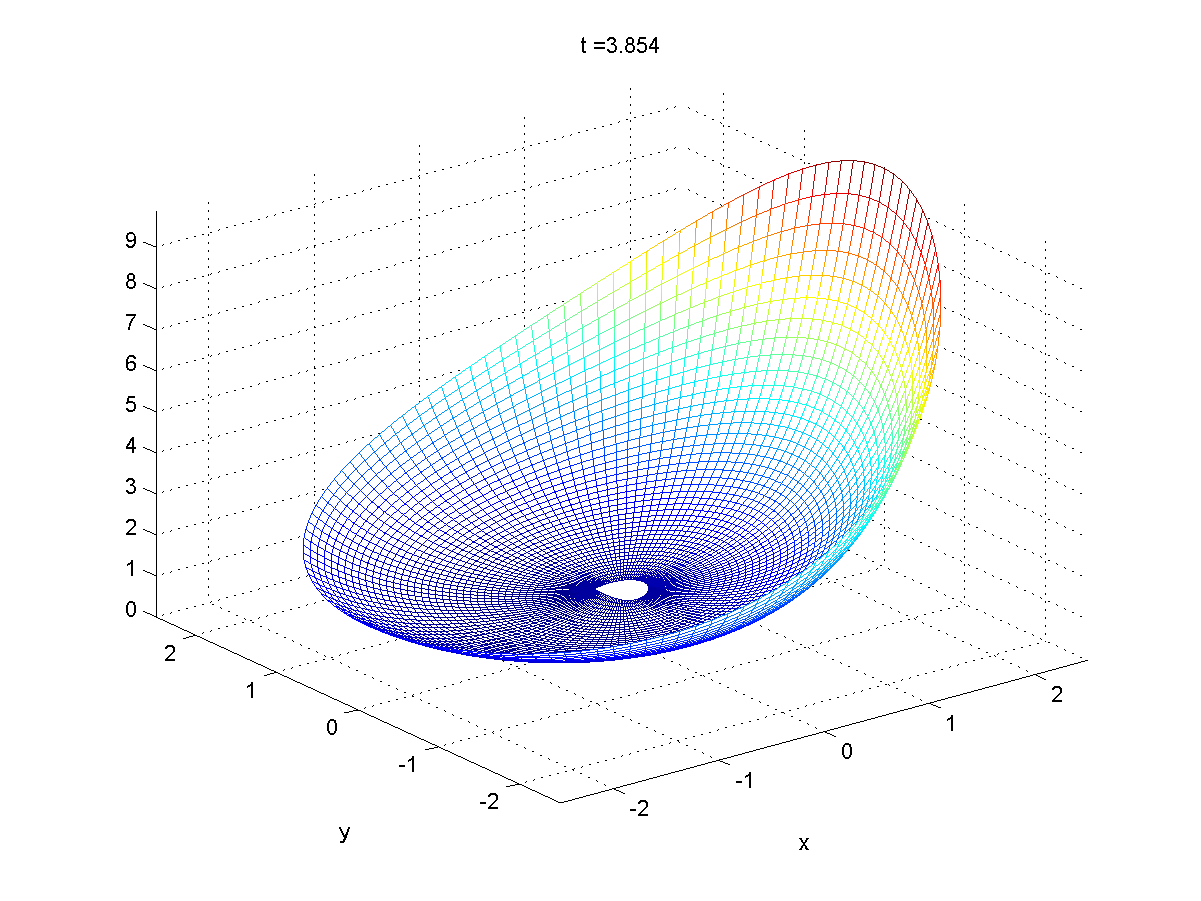} \\
   \includegraphics[scale=.25]{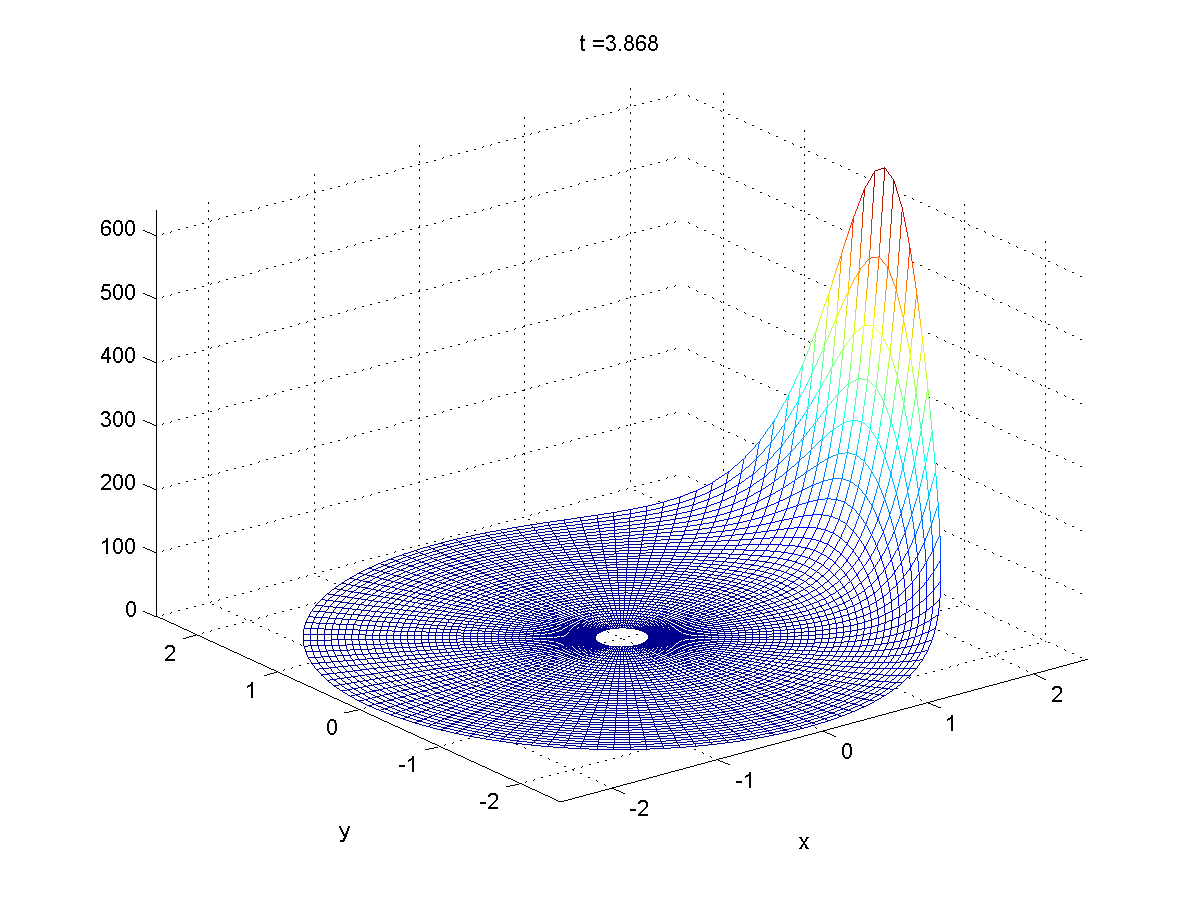}   & \includegraphics[scale=.25]{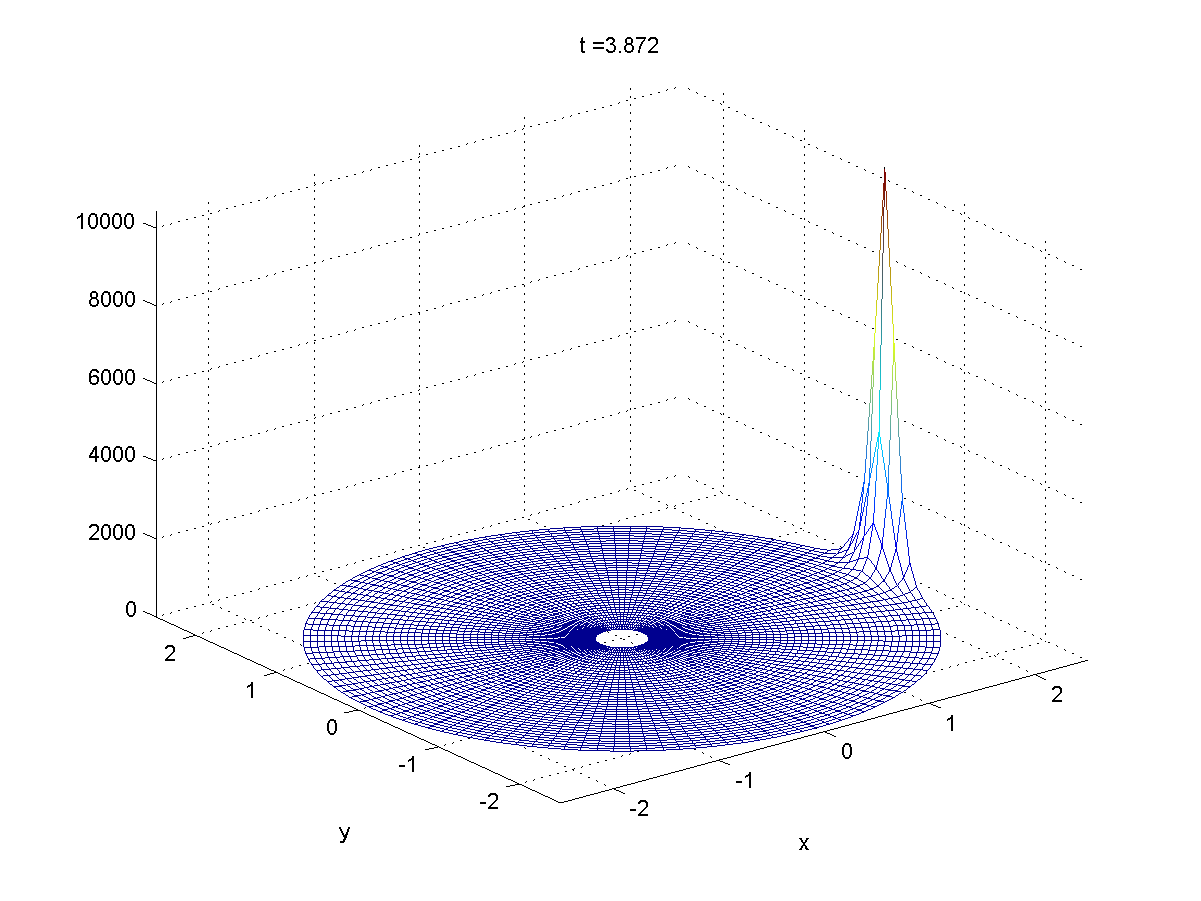}
\end{tabular}
 \caption{Numerical Simulations on $\Omega = [0.2,2.5] \times \mathbb{R}/2\pi \mathbb{Z}$ and all parameters equal to 1 with random initial conditions. For $M=10$ greater enough, symmetry breaking appears. Molecular markers are concentrated on one point of the membrane in finite time.}\label{figpolar1}
\end{figure}
\end{center}
The heuristic done in section \ref{Heuristic} allowed us thinking that the potential case could break symmetry more readily than the transversal case. We have also done numerical simulation in potential case for small $M$ and we found that the numerical behaviours are similar in the two possible drifts, see Fig \ref{figMT2} and \ref{figpolar2}.
We have assumed that the cell occupies a circle of radius $R>0$. Furthermore for simplicity, we consider a bounded-periodic domain $\Omega = [0,R] \times \mathbb{R}/2\pi R\mathbb{Z}$. With the numerical analysis done in \cite{proceedings} for the potential case, we see that the behaviours are similar with the annulus case developed in this article.

\section{Conclusion}

Polar description of the cell has been described in this article, this improvement fitted the real cell shape and was a first step before establishing a model for membrane deformation. In this work we have provided a first answer to the following question: do the nonlinear convection-diffusion models given in \cite{Firstpaper}  and \cite{Siam_CHMV} describe cell polarisation or not?  To do so we have used both a mathematical heuristic and numerical simulations, which have ensured us that solutions develop symmetry breaking over a critical value $M^*$. This has given us a first justification of the mathematical heuristic. On this point, the numerical behaviours are close to cell behaviours during biological experiences.  In order to fit biological measurements, the choice of parameters is essential and we refer to biological literature. Several measurements on polarisation time and localisation of polar cap have been made, we will describe them in a further work.\\

\noindent{\em Acknowledgement: this research has been supported by ANR program JCJC project MODPOL. N. Meunier and N. Muller want to thank M. de Buhan for helpful discussions.}

\medskip

\bibliographystyle{Siam}

\end{document}